# Binomial transform and the backward difference

Khristo N. Boyadzhiev
Ohio Northern University
Department of Mathematics and Statistics
Ada, Ohio 45810, USA
k-boyadzhiev@onu.edu

**Abstract** We show a remarkable property of the binomial transform – it converts multiplication by the discrete variable into a certain difference operator. We also consider the case of dividing by the discrete variable.
The properties presented here are used to compute various binomial transform formulas involving Harmonic numbers, Fibonacci numbers, Stirling numbers of the second kind, and Laguerre polynomials. Several new identities are proved and some known results are given new short proofs.

**1. Introduction and main results**

The two classical integral transforms, the Laplace and the Fourier transforms, have the important property that they convert, roughly speaking, multiplication by the variable into differentiation. In this paper we show that the discrete binomial transform has a similar property – it converts multiplication by the discrete variable $k$ into the operator $n\nabla$ where $\nabla$ is the backward difference. Because of the duality in the binomial transform, the converse is also true.

We also show what happens when we divide by the discrete variable.

Given a sequence $\{a_k\}_{k=0}^{\infty}$ we define its binomial transform to be the new sequence $\{b_n\}_{n=0}^{\infty}$ where

$$b_n = \sum_{k=0}^{n}\binom{n}{k}(-1)^{k-1}a_k \quad \text{with inversion} \quad \sum_{k=0}^{n}\binom{n}{k}(-1)^{k-1}b_k = a_n \ .$$

More information about the binomial transform can be found in [1], [5], and [10].

We want to see how $b_n$ changes when we multiply $a_k$ by $k$. Introducing the difference operator $\nabla b_n = b_n - b_{n-1}$ we have the following theorem.

**Theorem.** *For every positive integer $p$ and every $n \geq p$,*

$$\sum_{k=0}^{n} \binom{n}{k}(-1)^{k-1} k^p a_k = (n\nabla)^p b_n \ . \tag{1}$$

Here we need $n \geq p$ in the above formula because the RHS contains a term with $b_{n-p}$ which is not defined for $n < p$. We shall discuss this situation in the next section. It is clear that we only need to prove (1) for $p = 1$ and the rest follows by iteration.

We have the immediate corollary.

**Corollary 1.** *Let $g(t)$ be a polynomial. Then*

$$\sum_{k=0}^{n} \binom{n}{k}(-1)^{k-1} g(k) a_k = g(n\nabla) b_n \ .$$

In particular, for every complex number $\lambda$ and every integer $n \geq 1$,

$$\sum_{k=0}^{n} \binom{n}{k}(-1)^{k-1}(k+\lambda)\, a_k = (n\nabla + \lambda)b_n = (n+\lambda)b_n - nb_{n-1} \tag{2}$$

The presence of the factor $(-1)^{k-1}$ in (1) is unimportant. It is used here mainly because many popular formulas include this factor and the inversion formula with $(-1)^{k-1}$ is symmetrical. Because of this symmetry we can invert (1).

**Corollary 2.** *For every $n \geq 1$,*

$$\sum_{k=1}^{n} \binom{n}{k}(-1)^{k-1} k \nabla a_k = n b_n \ .$$

The next corollary describes division by $k + \lambda$.

**Corollary 3** *Let $a_0 = 0$ (so that $b_0 = 0$ too). Then for every number $\lambda \neq -1, -2, \ldots,$ and every $n \geq 1$,*

$$\sum_{k=1}^{n} \binom{n}{k}(-1)^{k-1} \frac{a_k}{k+\lambda} = \sum_{m=1}^{n} \frac{(m+1)(m+2)\ldots n}{(\lambda+m)(\lambda+m+1)\ldots(\lambda+n)} b_m \tag{3}$$

$$= n! \sum_{m=1}^{n} \frac{b_m}{m!(\lambda+m)(\lambda+m+1)\ldots(\lambda+n)} .$$

When $\lambda = 0$ equation (3) takes the form

$$\sum_{k=1}^{n} \binom{n}{k}(-1)^{k-1}\frac{a_k}{k} = \sum_{m=1}^{n} \frac{b_m}{m} . \qquad (4)$$

Property (4) was proved independently by A.N. 't Woord in [15]. When $a_0 = 0$, the two properties (1) and (4) are, in fact, equivalent – see Remark 1 at the end of the paper.

The case $\lambda = 1$ in (3) is also interesting,

$$\sum_{k=1}^{n} \binom{n}{k}(-1)^{k-1}\frac{a_k}{k+1} = \frac{1}{n+1}\sum_{m=1}^{n} b_m .$$

When $a_0 \neq 0, b_0 \neq 0$ we can apply (3) to the transform

$$\sum_{k=1}^{n} \binom{n}{k}(-1)^{k-1} a_k = b_n - b_0$$

and avoid the restriction. Thus we get, for instance (with summation on both sides from zero),

$$\sum_{k=0}^{n} \binom{n}{k}(-1)^{k-1}\frac{a_k}{k+1} = \frac{1}{n+1}\sum_{m=0}^{n} b_m . \qquad (5)$$

The above properties are true, of course, if we drop the alternating factor $(-1)^{k-1}$. That is, from the equation

$$\sum_{k=0}^{n} \binom{n}{k} a_k = b_n$$

it follows that

$$\sum_{k=0}^{n} \binom{n}{k} k^p a_k = (n\nabla)^p b_n$$

etc.

In the next section we shall present a number of examples. We want to demonstrate that properties (1), (3), (4), and (5) are very effective tools for evaluating binomial transforms and for generating various binomial transform formulas. The theorem and the last corollary will be proved in section 3.

## 2. Examples

Here we shall use the above results in order to obtain a number of new binomial transform identities and to give short proofs to some known ones. Our examples will involve harmonic, Stirling, and Fibonacci numbers, and also Laguerre polynomials.

The theorem can be used in two different ways: computing the binomial transform of $k^p a_k$ by computing $(n\nabla)^p b_n$, or computing $(n\nabla)^p b_n$ if the LHS in (1) can be evaluated by other means.

**Example 1**. We want to start with something simple. Consider the sequence $a_k = 1, k = 0,1,...$, where $b_0 = -1, b_n = 0, (n \geq 1)$, that is,

$$\sum_{k=0}^{n}\binom{n}{k}(-1)^{k-1} = 0 \ (n \geq 1) .$$

Then for $p=1$,

$$\sum_{k=0}^{n}\binom{n}{k}(-1)^{k-1}k = \begin{cases} 0 \ (n \neq 1) \\ 1 \ (n=1) \end{cases}$$

In general, for any positive integer $p$,

$$\sum_{k=0}^{n}\binom{n}{k}(-1)^{k-1}k^p = (-1)^{n-1}n!S(p,n). \tag{6}$$

where $S(p,n)$ are the Stirling numbers of the second kind with the property $S(p,n) = 0$ when $n > p$ and $S(n,n) = 1$. A very good reference for these numbers is the book [6]. Equation (6) is the classical representation of $S(p,n)$ as binomial transform.

**Example 2**. Let $a_k = x^k$, where $x$ is any real or complex number. Then

$$\sum_{k=0}^{n}\binom{n}{k}x^k = (1+x)^n$$

and therefore,

$$\sum_{k=0}^{n} \binom{n}{k} k^p x^k = (n\nabla)^p (1+x)^n .$$

At the same time

$$\sum_{k=0}^{n} \binom{n}{k} k^p x^k = \left(x\frac{d}{dx}\right)^p (1+x)^n = \sum_{j=0}^{p} S(p,j) x^j \left(\frac{d}{dx}\right)^j (1+x)^n \qquad (7)$$

$$= \sum_{j=0}^{p} S(p,j) x^j n(n-1)\ldots(n-j+1)(1+x)^{n-j} = \sum_{j=0}^{p} \binom{n}{j} S(p,j) j! x^j (1+x)^{n-j} .$$

We are using here the differentiation rule (see, for example, [10, p. 218] ).

$$\left(x\frac{d}{dx}\right)^p f(x) = \sum_{j=0}^{p} S(p,j) x^j \left(\frac{d}{dx}\right)^j f(x) .$$

The theorem then implies the following rule: For every for $n \geq p$,

$$(n\nabla)^p (1+x)^n = \sum_{j=0}^{p} \binom{n}{j} S(p,j) j! x^j (1+x)^{n-j} .$$

Now we shall involve property (4).

**Example 3.** Consider the sequence $a_0 = 0$, $a_k = 1 \, (n \geq 1)$ . Then $b_0 = 0$ and for $n \geq 1$

$$\sum_{k=1}^{n} \binom{n}{k} (-1)^{k-1} 1 = 1 .$$

From (4) we get

$$\sum_{k=1}^{n} \binom{n}{k} (-1)^{k-1} \frac{1}{k} = 1 + \frac{1}{2} + \ldots + \frac{1}{n} = H_n , \qquad (8)$$

where $H_n$ are the harmonic numbers, $H_0 = 0$. This formula is well-known. It appears, for instance, on p. 53 of Schwatt's book [11] and also on p. 6 in [4] and p. 5 in [10]. Property (4) makes it possible to give an immediate proof. Repeating (4) in (8) we also find,

$$\sum_{k=1}^{n} \binom{n}{k} (-1)^{k-1} \frac{1}{k^2} = \sum_{k=1}^{n} \frac{H_k}{k} ,$$

etc. Notice that by a simple computation

$$\sum_{k=1}^{n} \frac{H_k}{k} = \frac{1}{2}\left(H_n^2 + H_n^{(2)}\right), \text{ where } H_n^{(2)} = 1 + \frac{1}{2^2} + \ldots + \frac{1}{n^2}.$$

By inversion,

$$\sum_{k=1}^{n} \binom{n}{k}(-1)^{k-1}\left(H_k^2 + H_k^{(2)}\right) = \frac{2}{n^2}, \tag{9}$$

and again by (4)

$$\sum_{k=1}^{n} \binom{n}{k}(-1)^{k-1}\frac{1}{k}\left(H_k^2 + H_k^{(2)}\right) = 2\sum_{k=1}^{n}\frac{1}{k^3} = 2H_n^{(3)}.$$

Inverting this formula we get,

$$\sum_{k=1}^{n} \binom{n}{k}(-1)^{k-1} H_n^{(3)} = \frac{1}{2n}\left(H_n^2 + H_n^{(2)}\right).$$

Using now (1) in (9) we have for all $n > 1$

$$\sum_{k=1}^{n} \binom{n}{k}(-1)^{k-1} k \left(H_n^2 + H_n^{(2)}\right) = n\nabla \frac{2}{n^2} = \frac{2(1-2n)}{n(n-1)^2}, \tag{10}$$

$$\sum_{k=1}^{n} \binom{n}{k}(-1)^{k-1} k^2 \left(H_n^2 + H_n^{(2)}\right) = n\nabla\left(\frac{2(1-2n)}{n(n-1)^2}\right) = \frac{4}{(n-1)(n-2)}, \quad (n > 2),$$

$$\sum_{k=1}^{n} \binom{n}{k}(-1)^{k-1} k^3 \left(H_n^2 + H_n^{(2)}\right) = n\nabla\left(\frac{4}{(n-1)(n-2)}\right) = \frac{8n(n-2)}{(n-1)(n-3)}, \quad (n > 3),$$

etc. At the same time, inverting (8),

$$\sum_{k=1}^{n} \binom{n}{k}(-1)^{k-1} H_k = \frac{1}{n}$$

and applying (4)

$$\sum_{k=1}^{n} \binom{n}{k}(-1)^{k-1} \frac{H_k}{k} = H_n^{(2)}.$$

Inverting this we find

$$\sum_{k=1}^{n} \binom{n}{k}(-1)^{k-1} H_k^{(2)} = \frac{H_n}{n}, \tag{11}$$

and (1) yields consecutively,

$$\sum_{k=1}^{n}\binom{n}{k}(-1)^{k-1} k H_k^{(2)} = n\nabla\left(\frac{H_n}{n}\right) = \frac{1-H_n}{n-1} , \qquad (12)$$

$$\sum_{k=1}^{n}\binom{n}{k}(-1)^{k-1} k^2 H_k^{(2)} = n\nabla\left(\frac{1-H_n}{n-1}\right) = \frac{1-2n+nH_n}{(n-1)(n-2)} ,$$

etc. Subtracting (11) from (9) we find

$$\sum_{k=1}^{n}\binom{n}{k}(-1)^{k-1} H_k^2 = \frac{2}{n^2} - \frac{H_n}{n} ,$$

while subtracting (12) from (10) gives

$$\sum_{k=1}^{n}\binom{n}{k}(-1)^{k-1} k H_k^2 = \frac{H_n}{n-1} + \frac{2-3n-n^2}{n(n-1)^2} .$$

Proceeding further we can evaluate $\sum_{k=1}^{n}\binom{n}{k}(-1)^{k-1} k^p H_k^2$ for higher values of $p$. Chuanan Wei et al. have evaluated these sums up to $p=3$ in [14] by a different method.

**Example 4.** We continue to work here with the identity

$$\sum_{k=1}^{n}\binom{n}{k}(-1)^{k-1}\frac{1}{k^2} = \sum_{k=1}^{n}\frac{H_k}{k}$$

from the previous example. It can be written this way

$$\sum_{k=1}^{n}\binom{n}{k}(-1)^{k-1}\frac{1}{k^2} = \sum \frac{1}{k_1 k_2}, \quad 1 \le k_1 \le k_2 \le n .$$

Applying property (4) obviously yields

$$\sum_{k=1}^{n}\binom{n}{k}(-1)^{k-1}\frac{1}{k^3} = \sum \frac{1}{k_1 k_2 k_3}, \quad 1 \le k_1 \le k_2 \le k_3 \le n ,$$

and continuing this process we find for every positive integer $m$,

$$\sum_{k=1}^{n}\binom{n}{k}\frac{(-1)^{k-1}}{k^m} = \sum \frac{1}{k_1 k_2 ... k_m}, \quad 1 \le k_1 \le k_2 \le ... \le k_m \le n .$$

This formula was obtained by Karl Dilcher in 1995 by other means - see [1] for more information.

**Example 5**. Consider the sequence of harmonic numbers,

$$a_k = H_k = 1 + \frac{1}{2} + ... + \frac{1}{k} \quad (k > 0), \quad a_0 = 0, \text{ where } b_0 = 0 \text{ and for } n > 0,$$

$$\sum_{k=0}^{n} \binom{n}{k} (-1)^{k-1} H_k = \frac{1}{n} \quad .$$

Thus for $n > 1$,

$$\sum_{k=0}^{n} \binom{n}{k} (-1)^{k-1} k H_k = n\left(\frac{1}{n} - \frac{1}{n-1}\right) = \frac{-1}{n-1} \quad ,$$

and for $n > 2$,

$$\sum_{k=0}^{n} \binom{n}{k} (-1)^{k-1} k^2 H_k = n\left(\frac{-1}{n-1} + \frac{1}{n-2}\right) = \frac{n}{(n-1)(n-2)} \quad ,$$

etc. In fact, we have a general formula. From [1]

$$\sum_{k=1}^{n} \binom{n}{k} H_k k^p x^k = a(p,n,x) H_n - \sum_{k=1}^{n-1} \frac{a(p,k,x)}{n-k} \quad ,$$

where

$$a(p,n,x) = (x\frac{d}{dx})^p (x+1)^n = \sum_{j=0}^{n} \binom{n}{j} S(p,j) j! x^j (1+x)^{n-j} \quad .$$

With $x = -1$ we find $a(p,n,-1) = (-1)^n n! S(p,n)$ and thus

$$\sum_{k=1}^{n} \binom{n}{k} (-1)^{k-1} H_k k^p = (-1)^{n-1} n! S(p,n) H_n + \sum_{k=1}^{n-1} \frac{(-1)^k k! S(p,k)}{n-k} \quad . \tag{13}$$

The formula is true for any positive $p$ and $n$. When $n > p$ the first term on the RHS is missing, since $S(p,n) = 0$. In view of our theorem, this formula explains the action of the operator $(n\nabla)^p$ on the sequence $\{1/n\}_{n=1}^{\infty}$ for $n > p$, that is,

$$(n\nabla)^p \frac{1}{n} = \sum_{k=1}^{n-1} \frac{(-1)^k k! S(p,k)}{n-k} \quad .$$

**Example 6**. Starting from

$$\sum_{k=1}^{n}\binom{n}{k}(-1)^{k-1}x^k = 1-(1-x)^n$$

we find by inversion

$$\sum_{k=1}^{n}\binom{n}{k}(-1)^{k-1}\{1-(1-x)^k\} = x^n,$$

and then according to (4)

$$\sum_{k=1}^{n}\binom{n}{k}(-1)^{k-1}\left\{\frac{1-(1-x)^k}{k}\right\} = \sum_{k=1}^{n}\frac{x^k}{k}.$$

With $x = -1$ we get

$$\sum_{k=1}^{n}\binom{n}{k}(-1)^{k-1}\left\{\frac{1-2^k}{k}\right\} = \sum_{k=1}^{n}\frac{(-1)^k}{k} = -H_n^-,$$

where $H_n^- = 1 - \frac{1}{2} + \ldots + \frac{(-1)^{n-1}}{n}$ ($n \geq 1$) and $H_0^- = 0$ are the skew-harmonic numbers.

From here by inversion,

$$\sum_{k=1}^{n}\binom{n}{k}(-1)^{k-1}H_k^- = \frac{2^n-1}{n}. \tag{14}$$

According to our theorem

$$\sum_{k=1}^{n}\binom{n}{k}(-1)^{k-1}kH_k^- = n\nabla\left\{\frac{2^n-1}{n}\right\} = 2^{n-1}\frac{n-2}{n-1} + \frac{1}{n-1},$$

while (4) applied to (14) yields

$$\sum_{k=1}^{n}\binom{n}{k}(-1)^{k-1}\frac{H_k^-}{k} = \sum_{k=1}^{n}\frac{2^k}{k^2} - H_n^{(2)}.$$

With (5) applied to (14),

$$\sum_{k=1}^{n}\binom{n}{k}(-1)^{k-1}\frac{H_k^-}{k+1} = \frac{1}{n+1}\left\{\sum_{k=1}^{n}\frac{2^k}{k} - H_n\right\}.$$

**Example 7**. The operator $n\nabla$ works very well on sequences defined by recurrence relations. For example, let $F_0, F_1, F_2,...$ be the sequence of Fibonacci numbers, where $F_n = F_{n-1} + F_{n-2}$ (the books [7] and [13] are good references). We shall prove some new identities. It is known that for every $n \geq 0$,

$$\sum_{k=0}^{n}\binom{n}{k}(-1)^{k-1}F_k = F_n$$

(see [5] ). Thus from (1),

$$\sum_{k=0}^{n}\binom{n}{k}(-1)^{k-1}kF_k = n(F_n - F_{n-1}) = nF_{n-2}$$

and further,

$$\sum_{k=0}^{n}\binom{n}{k}(-1)^{k-1}k^2 F_k = n^2 F_{n-4} + nF_{n-3} ,$$

etc. At the same time, without the alternating factor $(-1)^{k-1}$ we have the classical result of Edouard Lucas [7],

$$\sum_{k=0}^{n}\binom{n}{k}F_k = F_{2n} .$$

From here,

$$\sum_{k=0}^{n}\binom{n}{k}kF_k = n(F_{2n} - F_{2n-2}) = nF_{2n-1} ,$$

$$\sum_{k=0}^{n}\binom{n}{k}k^2 F_k = n^2 F_{2n-2} + nF_{2n-3} ,$$

etc. Almost the same equations are satisfied by the Lucas numbers $L_n$, $n \geq 0$ ( $2,1,3,4,...$ ), as they are defined by the same recurrence relation $L_n = L_{n-1} + L_{n-2}$ and

$$\sum_{k=0}^{n}\binom{n}{k}(-1)^k L_k = L_n , \quad \sum_{k=0}^{n}\binom{n}{k}L_k = L_{2n} .$$

For the Fibonacci numbers we also have from (5)

$$\sum_{k=0}^{n}\binom{n}{k}(-1)^{k-1}\frac{F_k}{k+1}=\frac{1}{n+1}\sum_{k=0}^{n}F_k=\frac{1}{n+1}(F_{n+2}-1),$$

by using the fundamental property $F_1+F_2+...+F_n=F_{n+2}-1$ and the fact that $F_0=0$. For the Lucas numbers we also have $L_0+L_1+L_2+...+L_n=L_{n+2}-1$ and therefore,

$$\sum_{k=0}^{n}\binom{n}{k}(-1)^k\frac{L_k}{k+1}=\frac{1}{n+1}(L_{n+2}-1).$$

**Example 8**. We shall use now a recurrence property to evaluate one special binomial transform. Let $q$ be a nonnegative integer and set for $n\geq 1$

$$\sigma_n(q)=1^q+2^q+...+n^q,\ \sigma_0(q)=0.$$

We want to compute the sequence

$$b_n=\sum_{k=0}^{n}\binom{n}{k}(-1)^{k-1}\sigma_k(q),\ (n=0,1,...).$$

Obviously, $b_0=0$ and $b_1=1$. Next we invert the above equation to get

$$\sum_{k=0}^{n}\binom{n}{k}(-1)^{k-1}b_k=\sigma_n(q),$$

and use property (1) for $n\geq 1$

$$\sum_{k=0}^{n}\binom{n}{k}(-1)^{k-1}kb_k=n(\sigma_n(q)-\sigma_{n-1}(q))=n^{q+1}.$$

At the same time, by inversion in (6),

$$\sum_{k=0}^{n}\binom{n}{k}k!S(q+1,k)=n^{q+1}. \qquad (15)$$

Therefore, $(-1)^{k-1}kb_k=k!S(q+1,k)$ and $b_k=(-1)^{k-1}(k-1)!S(q+1,k)$, $k=1,2,....$ . That is,

$$\sum_{k=0}^{n}\binom{n}{k}(-1)^{k-1}\sigma_k(q)=(-1)^{n-1}(n-1)!S(q+1,n)$$

for $n\geq 1$. It is remarkable that this sequence truncates, $b_n=0$ for $n>q+1$.

**Example 9**. This example is related to the previous example. For $q \geq 1$ we can write (15) in the form

$$\sum_{k=1}^{n}\binom{n}{k}k!S(q,k) = n^q ,$$

since $S(q,0) = 0$ for $q \geq 1$. Then property (5) provides the new formula.

$$\sum_{k=1}^{n}\binom{n}{k}\frac{k!S(q,k)}{k+1} = \frac{1}{n+1}\sum_{k=1}^{n}k^q = \frac{\sigma_n(q)}{n+1} .$$

**Example 10.** We shall use again the identity

$$\sum_{k=1}^{n}\binom{n}{k}(-1)^{k-1}H_k = \frac{1}{n} \quad (n \geq 1)$$

(the inverse of (8)) and apply property (3). For any $\lambda \neq -1, -2, ...,$ we have

$$\sum_{k=1}^{n}\binom{n}{k}(-1)^{k-1}\frac{H_k}{k+\lambda} = \sum_{m=1}^{n}\frac{(m+1)(m+2)...n}{(\lambda+m)(\lambda+m+1)...(\lambda+n)m} . \tag{16}$$

For $\lambda = 1$ this turns into the symmetric formula

$$\sum_{k=1}^{n}\binom{n}{k}(-1)^{k-1}\frac{H_k}{k+1} = \frac{H_n}{n+1} ,$$

which is an example of a sequence invariant under the binomial transform.

For $\lambda = 2$, $\lambda = 3$, and $\lambda = 4$ in (16) we have correspondingly.

$$\sum_{k=1}^{n}\binom{n}{k}(-1)^{k-1}\frac{H_k}{k+2} = \frac{H_n + n}{(n+1)(n+2)} ,$$

$$\sum_{k=1}^{n}\binom{n}{k}(-1)^{k-1}\frac{H_k}{k+3} = \frac{n^2 + 7n + 4H_n}{2(n+1)(n+2)(n+3)} ,$$

$$\sum_{k=1}^{n}\binom{n}{k}(-1)^{k-1}\frac{H_k}{k+4} = \frac{2n^3 + 21n^2 + 85n + 36H_n}{6(n+1)(n+2)(n+3)(n+4)} .$$

A recurrence relation for the evaluation of the sums in (16) was developed in [3] and the authors have computed these sums for $\lambda = 1, 2, 3, 4$. The last two results (for $\lambda = 3, 4$) in [3, p. 2227], however, are wrong; they contain additional terms in the numerators.

**Example 11.** Consider the sequence $a_k = k$, $k = 1, 2, ...$, where $b_1 = 1$ and $b_n = 0$ for $n > 1$ (see Example 1). Suppose $\lambda \neq 0$. We apply property (3) to the equation

$$\sum_{k=1}^{n} \binom{n}{k}(-1)^{k-1} k = \begin{cases} 0 \ (n > 1) \\ 1 \ (n = 1) \end{cases}$$

to get

$$\sum_{k=1}^{n} \binom{n}{k}(-1)^{k-1} \frac{k}{k+\lambda} = \frac{n!}{(\lambda+1)(\lambda+2)...(\lambda+n)} \ .$$

From here, writing $\dfrac{k}{k+\lambda} = \dfrac{k+\lambda-\lambda}{k+\lambda} = 1 - \dfrac{\lambda}{k+\lambda}$ we find

$$\frac{1}{\lambda} + \sum_{k=1}^{n} \binom{n}{k}(-1)^{k} \frac{1}{k+\lambda} = \frac{n!}{\lambda(\lambda+1)(\lambda+2)...(\lambda+n)} \ ,$$

or, starting he summation from $k = 0$,

$$\sum_{k=0}^{n} \binom{n}{k}(-1)^{k} \frac{1}{k+\lambda} = \frac{n!}{\lambda(\lambda+1)(\lambda+2)...(\lambda+n)} \ ,$$

which is a well known important identity [4], [6, p.188].

**Example 12.** Let $L_n(x) = \dfrac{e^x}{n!}\left(\dfrac{d}{dx}\right)^n (x^n e^{-x})$ be the Laguerre polynomials. It is known that

$$\sum_{k=0}^{n} \binom{n}{k} \frac{(-x)^k}{k!} = L_n(x)$$

(see p. 213 in [9]), which can be written also in the form (starting summation from $k = 1$)

$$\sum_{k=1}^{n} \binom{n}{k} \frac{(-x)^k}{k!} = L_n(x) - 1. \tag{17}$$

Using property (4) in (17) we find a new identity

$$\sum_{k=1}^{n} \binom{n}{k} \frac{(-x)^k}{k!\,k} = \sum_{k=1}^{n} \frac{L_k(x)}{k} - H_n \ .$$

The LHS here can be obtained also by dividing (17) by $x$ and integrating. Thus

$$\sum_{k=1}^{n}\binom{n}{k}\frac{(-x)^k}{k!k}=\int_{0}^{x}\frac{L_k(t)-1}{t}dt \ .$$

Therefore (cf. [2]),

$$\sum_{k=1}^{n}\frac{L_k(x)}{k}=\int_{0}^{x}\frac{L_k(t)-1}{t}dt+H_n \ .$$

Next, applying (5) to (17) yields

$$\sum_{k=1}^{n}\binom{n}{k}\frac{(-x)^k}{k!(k+1)}=\frac{1}{n+1}\sum_{k=1}^{n}(L_k(x)-1) \ .$$

and then by integration in (17) we obtain the property

$$\frac{1}{n+1}\sum_{k=1}^{n}(L_k(x)-1)=\frac{1}{x}\int_{0}^{x}L_n(t)-1\,dt \ .$$

**Example 13**. Let $p$ be a nonnegative integer. Consider the binomial transform

$$\sum_{k=0}^{n}\binom{n}{k}\binom{p}{k}=\binom{p+n}{p},$$

which is known as the Vandermonde identity [4], [6]. According to (1),

$$\sum_{k=0}^{n}\binom{n}{k}\binom{p}{k}k=n\left\{\binom{p+n}{p}-\binom{p+n-1}{p}\right\}=n\binom{p+n-1}{p-1}.$$

Applying (5) we find

$$\sum_{k=0}^{n}\binom{n}{k}\binom{p}{k}\frac{1}{k+1}=\frac{1}{n+1}\sum_{j=0}^{n}\binom{p+j}{p}=\frac{1}{n+1}\binom{p+n+1}{p+1} \ .$$

### 3. Proofs

First we shall prove the theorem. For simplicity we drop here the factor $(-1)^{k-1}$. Thus we only need to show that if for every $n \geq 0$ we have

$$\sum_{k=0}^{n}\binom{n}{k}a_k=b_n \ ,$$

then for $n \geq 1$,

$$\sum_{k=0}^{n} \binom{n}{k} k a_k = n(b_n - b_{n-1}) \ .$$

First, notice that for any $n \geq 1$ and any $0 \leq k \leq n$,

$$\binom{n}{k} - \binom{n-1}{k} = \frac{k}{n} \binom{n}{k}$$

Therefore, we have

$$b_n - b_{n-1} = \sum_{k=0}^{n} \left[ \binom{n}{k} - \binom{n-1}{k} \right] a_k = \sum_{k=0}^{n} \left[ \frac{k}{n} \binom{n}{k} \right] a_k = \frac{1}{n} \sum_{k=0}^{n} \binom{n}{k} k a_k$$

and the assertion follows for $p = 1$. Repeating this $p > 1$ times we obtain (1).

*Proof of Corollary 3.*

Take an arbitrary $\lambda$ and define the sequence $\{c_n\}_{n=1}^{\infty}$ by the equation

$$\sum_{k=1}^{n} \binom{n}{k} (-1)^{k-1} \frac{a_k}{k+\lambda} = c_n \ .$$

When $n=1$ we see directly from here that $c_1 = \dfrac{b_1}{\lambda+1}$. For $n > 1$ we find from (2)

$$\sum_{k=1}^{n} \binom{n}{k} (-1)^{k-1} a_k = (n+\lambda) c_n - n c_{n-1} = b_n \ ,$$

and so we can compute $c_n$ from the recurrence relation

$$(n+\lambda) c_n - n c_{n-1} = b_n$$

for $n \geq 2$. Thus for $n = 2$ we have $(\lambda + 2) c_2 - 2 c_1 = b_2$ and therefore,

$$c_2 = \frac{2 b_1}{(\lambda+1)(\lambda+2)} + \frac{b_2}{\lambda+2} \ .$$

For $n = 3$ we find in the same manner

$$c_3 = \frac{2\cdot 3}{(\lambda+1)(\lambda+2)(\lambda+3)}b_1 + \frac{3}{(\lambda+2)(\lambda+3)}b_2 + \frac{1}{\lambda+3}b_3,$$

and (3) follows by induction. The details are left to the reader.

**Remark 1.** When $a_0 = 0$ (and hence $b_0 = 0$) properties (1) and (4) are equivalent. We just need to see that (1) follows from (4). Here we shall sketch the proof. Suppose (4) is true. Then let for $n = 1, 2, ...,$ we define the sequence $c_n$ by

$$\sum_{k=1}^{n}\binom{n}{k}(-1)^{k-1} k\, a_k = c_n.$$

Clearly, $c_1 = b_1$. For $n > 1$ from (4),

$$\sum_{k=1}^{n}\binom{n}{k}(-1)^{k-1} a_k = \sum_{m=1}^{n}\frac{c_m}{m} = b_n,$$

and now we can compute $c_n$ in terms of $b_n$ and $b_{n-1}$ from the recurrence relation

$$\sum_{m=1}^{n}\frac{c_m}{m} = b_n. \tag{18}$$

For $n = 2$ we have $\frac{c_2}{2} + c_1 = b_2$, that is, $c_2 = 2(b_2 - b_1)$. It is easy to see that if (18) holds and $c_m = m(b_m - b_{m-1})$ is true for $m \leq n$, then it is also true also for $m = n+1$. The simple algebra is again left to the reader.

**Remark 2**. Interesting results about binomial formulas involving the forward difference operator $\Delta$ were obtained by Spivey [12].

**Remark 3**. The backward difference operator $\nabla$ was used by Nielsen in [8] to study Bernoulli numbers. It was discussed also in [10]. The author hopes that the above results will bring new life to this operator.